\documentclass[referee,sn-mathphys]{sn-jnl_a}

\usepackage{amsthm,amsmath,amsfonts,amssymb}

\usepackage{mathtools}
\usepackage{blkarray}
\usepackage{xcolor}
\usepackage{subfig}
\usepackage{colortbl}

\usepackage{algorithm,algpseudocode}

\newcommand*{\Scale}[2][4]{\scalebox{#1}{$#2$}}%

\DeclareMathOperator{\sign}{sgn}

\jyear{2021}%

\theoremstyle{thmstyleone}%
\newtheorem{theorem}{Theorem}[section]
\newtheorem{lemma}[theorem]{Lemma}

\theoremstyle{definition}

\newtheorem{example}[theorem]{Example}

\theoremstyle{remark}

\numberwithin{equation}{section}

\begin{document}


\title[An Algorithm for Approximating Implicit Functions by Polynomials]{An Algorithm for Approximating Implicit Functions by Polynomials without Higher-Order Differentiability}


\author*[1]{\fnm{Kyung Soo} \sur{Rim}}\email{ksrim@sogang.ac.kr}



\affil*[1]{\orgdiv{Department of Mathematics}, \orgname{Sogang University}, \orgaddress{\street{35 Baekbeom-ro, Mapo-gu}, \city{Seoul}, \postcode{04107}, \country{Korea}}}

%


\abstract{%
We consider an equation of multiple variables in which a partial derivative does not vanish at a point. The implicit function theorem provides a local existence and uniqueness of the function for the equation. In this paper, we propose an algorithm to approximate the function by a polynomial without using higher-order differentiability, which depends essentially on integrability. Moreover, we extend the method to a system of equations if the Jacobian determinant does not vanish. This is a robust method for implicit functions that are not differentiable to higher-order. Additionally, we present two numerical experiments to verify the theoretical results.
}




\maketitle

\section{Introduction} 
Implicitly induced equations have been studied extensively because many useful mathematical models have expressions in the form of equations of several variables. A central issue in this subject is how to separate dependent variables as functions from the equations. This variable separation enables us to estimate mathematical models more easily. Specifically, for a continuously differentiable function 
$f:\mathbb{R}^2\times\mathbb{R}\to\mathbb{R}$, the problem is to solve 
\begin{equation} \label{main-problem}
	f(x,y)=0
\end{equation}
for $y$ which is a function of $x$ in a rectangle.

The problem of finding an implicit function defined by (\ref{main-problem}) has a long history dating back to the 1660s, and has been studied extensively (\cite{newton}, \cite{struik}, \cite[Theorem 2.3.1]{kp}). Various aspects of this problem, such as analyzing the behavior or justifying the existence of an implicitly defined function, have been rigorously formulated and studied since the 19th century (\cite[Theorems 2.4.6 and 6.1.2]{kp}, \cite{dini}).

The development of the implicit function theory has yielded many particular types of results that have found application in various fields, including mechanics, physics, engineering, economics, and mathematics. The theory has also been extended to Banach spaces, even under degenerate or non-smooth situations, based on a profound mathematical background (refer to, e.g., \cite{nash,moser,kp,dr}).

However, the focus of this study is to approximate an explicit form of $y=g(x)$ such that (\ref{main-problem}) can be solved without requiring deep mathematical knowledge and is easy to implement for practical applications. In Section 5, we demonstrate the computational validity of the proposed method.

In this paper, absolute values appear for certain quantities, where if the quantity is a scalar, it represents the classical absolute value, if it is a matrix, it represents the determinant, and if it is a subset of Euclidean space, it represents length, area, or volume.

\section{Integral of Heaviside composition} \label{implicit functions}

As a preliminary step, we will integrate the composition of a Heaviside function with $f$ in (\ref{main-problem}). We assume the existence of a point $(a,b)$ such that $\partial_yf(a,b)\ne0$ and $f(a,b)=0$. By applying the implicit function theorem, we obtain a rectangle $U$ containing $a$, an interval $V$ containing $b$, and a unique continuously differentiable function $g:U\to V$ such that $f(x,g(x))=0$ holds for all $x$ in $U$. Note that all (sub)rectangles in this paper are parallel to the coordinate axes. Therefore, we will restrict the domain and range of $f$ to $U\times V\to V$ going forward.

Let $U$ and $V$ be fixed. For any subrectangle $R$ of $U$, we define the integral of a composite function as follows:
\begin{equation} \label{volume integration in poly}
\mu(R)=\int_{R\times V} \Theta(f(x,y)) dxdy,
\end{equation}
where, $\Theta$ is the Heaviside function, which is defined as $\Theta(t)=1$ if $t\ge0$ and $\Theta(t)=0$ otherwise. The integral in (\ref{volume integration in poly}) is taken with respect to the standard Lebesgue measure $dxdy$. For a more detailed implementation of (\ref{volume integration in poly}), see \cite{wen}. Note that the quantity in (\ref{volume integration in poly}) will serve as input data to identify the function $y=g(x)$.

Fix $x\in U$. Since the sign function $\mathrm{sign}(f(x,y))$ is a single-step function on $V$, we can define $\rho=\rho(x)$ to be $1$ if $\mathrm{sign}(f(x,y))$ is increasing and $-1$ if it is decreasing. Then, from the continuity of $f$, we see that $\rho$ is a constant on $U\times V$, taking values of either $1$ or $-1$.

We now state the following lemma, which relates (\ref{volume integration in poly}), (\ref{main-problem}), and the function $y=g(x)$:

\begin{lemma} \label{coro of representation}
Under the above notations of $f$, $U$, $V$, $\mu$, and $\rho$, let $g:U\to V$ be such that $f(x,g(x))=0$ for all $x\in U$.  Then for any subrectangle $R$ of $U$,
\begin{equation*}
	\int_R gdx=\vert R\vert (\max V)^{\frac{\rho+1}{2}}(\min V)^{\frac{1-\rho}{2}} - \rho\mu(R).
\end{equation*}
\end{lemma}

\begin{proof}
Let $R$ be a subrectangle of $U$. First, suppose $\rho=1$. Since $\mu(R)$ is equal to the volume of the region enclosed by the planes $y=\max V$ and $y=g(x)$ in $R\times V$, we have
\begin{equation} \label{sign 1}
	\mu(R)= \vert R\vert \max V - \int_R gdx.
\end{equation}

If $\rho=-1$, then similarly,
\begin{equation} \label{sign 1-1}
	\mu(R) = \vert R\vert (-\min V) + \int_R gdx.
\end{equation}

Combining (\ref{sign 1}) and (\ref{sign 1-1}), and multiplying $\rho$ on both sides, we get the desired formula.
\end{proof}

\section{Step function by local averages} \label{aif}
In this section, we will construct a sequence of polynomials such that their local averages converge to $y=g(x)$.
We begin by introducing some terminology about partitions of $U$. Let $\Gamma_n$ be a partition of $U$ such that $\Gamma_{n+1}$ is finer than $\Gamma_n$ for every nonnegative integer $n$, and let $\Gamma = \bigcup_n \Gamma_n$. We say that $\Gamma$ is \emph{shrinkable} on $U$ if, for any point $x$ in $U$, there exists a sequence $R_n$ of blocks from $\Gamma_n$ that shrink to $x$ as $n \to \infty$. We refer to such a sequence $R_n$ as an \emph{$x$-shrinking block}.

Suppose that $\Gamma$ is shrinkable. For a sequence of integrable functions $g_n$, we say that $g_n$ converges to $g$ in the \emph{$\Gamma$-Cesàro sense} if, for any $x \in U$ and any $\epsilon > 0$, there exists an $N$ such that, for every $x$-shrinking block $R_n$, we have
\begin{equation} \label{spatial cesaro}
	\left\vert\tilde g_n(x)-g(x)\right\vert<\epsilon
\end{equation}
for all $n \geq N$, where $\tilde g_n(x) = \frac{1}{\vert R_n\vert} \int_{R_n} g_n  dx$ for $x \in R_n$. Note that $\tilde g_n$ is a step function that is constant on each block of $\Gamma_n$. For more information on the average operators of $\tilde g_n$, refer to \cite{jmz, folland, grafa, stein, zygmund, hk}.

To implement, we construct a shrinkable collection by dyadic decomposition. For notational simplicity, suppose that $U$ is a unit cube with center $a$, i.e., $U=[0,1]^2+a$, where $(a,b)$ is as in Section \ref{implicit functions} and satisfies $\partial_y f(a,b)\ne0$ with $f(a,b)=0$. Let $R_{n,i}=[(i-1)/2^n,i/2^n]\subset[0,1]$ $(1\le i\le 2^n)$. We set $\Gamma_n=\{R_{n,(i,j)}=R_{n,i}\times R_{n,j}\subset [0,1]^2\,:\,1\le i,j\le 2^n\}$. Then $\Gamma$ is shrinkable on $U$. For $a=(a_1,a_2)$, we define auxiliary matrices $A_n(a_k)$ whose $(i,j)$-component is the integral of
\begin{equation} \label{alpha}
 	\int_{R_{n,i}} (t-a_k)^j dt \qquad (1\le i, j\le 2^n,\;k=1,2)
\end{equation}
and $B_n$ whose $(i,j)$-component is calculated as
\begin{equation} \label{beta}
 	\vert R_{n,(i,j)}\vert (\max V)^{\frac{\rho+1}{2}}(\min V)^{\frac{1-\rho}{2}}-\rho\mu(R_{n,(i,j)}) \qquad(1\le i, j\le 2^n),
\end{equation}
where $\mu(R_{n,(i,j)})$ follows from Lemma \ref{coro of representation}.

For real numbers $c_{\alpha,\beta}$, putting
\begin{equation} \label{seq poly}
	g_n(x)=\sum_{\alpha=0}^{2^n-1}\sum_{\beta=0}^{2^n-1}c_{\alpha,\beta}(x_1-a_1)^\alpha(x_2-a_2)^\beta,
\end{equation}
we state and prove the main result of this section. 

\begin{theorem} \label{prob main1}
If the coefficient matrix $(c_{\alpha,\beta})$ of (\ref{seq poly}) is set to be 
\begin{equation*}
	\rho A_n(a_1)^{-1}B_n (A_n(a_2)^{T})^{-1},
\end{equation*}
then $\tilde g_{n}$ converges to $g$ as $n\to\infty$.
\end{theorem}

\begin{proof} 
We will prove the statement in three steps. For a positive integer $n$ let $\Delta=2^{-n}$.

Step 1: We will show that $A_n(a_k)$ is invertible and that its determinant is independent of $a_k$. By directly evaluating (\ref{alpha}) and using the invariant property of determinants, we obtain the following determinant of the Vandermonde matrix:
\begin{equation} \label{determinant}
\begin{split}
\vert A_n(a_k)\vert &=
	\arraycolsep=3pt\def\arraystretch{1.5}
	\frac{1}{2^n!}\left\vert\;
	\begin{matrix} 
	\Delta & \Delta^2  & \cdots & \Delta^{2^n} \\
	2\Delta & (2\Delta)^2 & \cdots & (2\Delta)^{2^n} \\
	\vdots & \vdots & \ddots & \vdots  \\
	 2^n\Delta & (2^n\Delta)^2 & \cdots & (2^n\Delta)^{2^n} \\
	\end{matrix}\right\vert \\  
	&=
	\arraycolsep=6pt\def\arraystretch{1.5}
	\Delta^{2^n}\left\vert\;
	\begin{matrix} 
	1 & \Delta  & \cdots & \Delta^{2^n-1} \\
	1 & 2\Delta & \cdots & (2\Delta)^{2^n-1} \\
	\vdots & \vdots & \ddots & \vdots  \\
	 1 & 2^n\Delta & \cdots  & (2^n\Delta)^{2^n-1} \\
	\end{matrix}\right\vert \\
	&=
 	\Delta^{2^n}\prod_{1\le i<j\le 2^n} (j-i)\Delta,
\end{split}	
\end{equation}
where the last term is strictly positive, which is independent of $a_k$.

\smallskip

Step 2. We will obtain a unique coefficients matrix $(c_{\alpha,\beta})$ of $g_n$ such that
\begin{equation}\label{gamma}
	\int_{R_{n,(i,j)}}g_ndx = \int_{R_{n,(i,j)}}gdx
\end{equation} 
for every $1\le i,j\le 2^n$. By Lemma \ref{coro of representation},
\begin{equation*} \label{bridge1}
	\int_{R_{n,(i,j)}}g dx = 
	\vert R_{n,(i,j)}\vert (\max V)^{\frac{\rho+1}{2}}(\min V)^{\frac{1-\rho}{2}} - \rho\mu(R_{n,(i,j)}).
\end{equation*} 
Using (\ref{gamma}) and the above equation, we get
\begin{equation} \label{bridge}
	\int_{R_{n,(i,j)}}g_n dx = 
	\vert R_{n,(i,j)}\vert (\max V)^{\frac{\rho+1}{2}}(\min V)^{\frac{1-\rho}{2}} - \rho\mu(R_{n,(i,j)}).
\end{equation} 
Substituting the definition of $A_n(a_k)$ into (\ref{bridge}), we obtain
\begin{equation*} 
\sum_{\alpha=0}^{2^n-1}\sum_{\beta=0}^{2^n-1}c_{\alpha,\beta} \frac{v_{i,\alpha}(a_1)}{\alpha+1}\frac{v_{j,\beta}(a_2)}{\beta+1}
	= \Delta^2(\max V)^{\frac{\rho+1}{2}}(\min V)^{\frac{1-\rho}{2}} - \rho\mu(R_{n,(i,j)}),
\end{equation*}
where $v_{i,j}(a_k)$ is the $(i,j)$-component of $A_n(a_k)$.
This can be expressed as a matrix equation, i.e.,
\begin{equation*}
	A_n(a_1)(c_{\alpha,\beta})A_n(a_2)^T = B_n.
\end{equation*}
where $B_n$ is the matrix whose $(i,j)$-entry is given by
\begin{equation*}
	\Delta^2(\max V)^{\frac{\rho+1}{2}}(\min V)^{\frac{1-\rho}{2}} - \rho\mu(R_{n,(i,j)}).
\end{equation*}
Since $A_n(a_k)$ is invertible by Step 1, we can solve for $(c_{\alpha,\beta})$ as follows:
\begin{equation*}
(c_{\alpha,\beta}) = A_n(a_1)^{-1}B_n(A_n(a_2)^{T})^{-1}.
\end{equation*}

\smallskip

Step 3. We will show that $g_n$ induced by Step 2 converges to $g$ in the $\Gamma$-Ces\`aro sense.
Let $x_0\in U$. Since $\Gamma$ is shrinkable, there exists a collection of rectangles $R_n\in\Gamma$ that shrinks to $x_0$. Since $g$ is continuous, by the Lebesgue differentiation theorem, for any $\epsilon>0$, there exists an $N$ such that for every $n\ge N$,
\begin{equation} \label{mean}
	\left\vert\frac{1}{\vert R_n\vert}\int_{R_n}gdx -g(x_0)\right\vert<\epsilon
\end{equation}

For $n\ge N$, take $m$ such that $m\ge n$ so that there is a collection of non-overlapping rectangles $R_{m,(i,j)}\in\Gamma_m$ such that
\begin{equation*}
	R_n= \bigcup_{(i,j)}R_{m,(i,j)}.
\end{equation*}
It follows that
\begin{equation*}
\begin{split}
	\int_{R_n} gdx 
	&= \sum_{(i,j)}\int_{R_{m,(i,j)}}gdx \\
	&= \sum_{(i,j)}\int_{R_{m,(i,j)}}g_mdx  \qquad \mbox{by (\ref{gamma}) of Step 2} \\	
	&= \int_{R_{n}}g_mdx
\end{split}
\end{equation*}
for every $m\ge n\ge N$. Taking $m=n$ and combining this with (\ref{mean}), we conclude that $g_n$ converges to $g$ in the $\Gamma$-Cesàro sense. Therefore, the proof is complete.
\end{proof}

If $g$ is a polynomial in Theorem \ref{prob main1}, then the sequence $g_n$ defined by Step 2 of the proof coincides with $g$ for all $n$ greater than or equal to the largest exponent of each variable in $g$.

\section{System of implicit functions} \label{section of system idv}

As a generalization of Theorem \ref{prob main1} in Section \ref{aif}, we consider the system of two equations given by
\begin{equation} \label{system idrv}
	f(x,y)=(f_1(x,y),f_2(x,y))=(0,0)
\end{equation}
with $f(a,b)=(0,0)$, where $f:\mathbb{R}^2\times\mathbb{R}^2\to\mathbb{R}^2$.
We aim to approximate the vector-valued function $y=(p(x),q(x))$ such that $f(x,p(x),q(x))=(0,0)$ in some neighborhood of $(a,b)$, where $p,q:\mathbb{R}^2\to\mathbb{R}$.

Assuming that $f$ is continuously differentiable and there exists a point $(a,b)$ such that $\det J_{f,y}(a,b)\neq0$ and $f(a,b)=(0,0)$, where $J_{f,y}(a,b)$ is the Jacobian matrix of $f$ with respect to $y$ at $(a,b)$, we can apply the implicit function theorem. This theorem guarantees that there is a unique vector-valued function $y=(p(x),q(x)) : U \to V$, defined on rectangles $U$ containing $a$ and $V$ containing $b$, such that $f(x,p(x),q(x))=(0,0)$ in $U$ and $b=(p(a),q(a))$.

First, we will state and prove the following lemma.

\begin{lemma} \label{system lemma}
Suppose that $\det J_{f,y}(a,b)\neq0$ with the same notation as above. Then:
\begin{itemize}
\item[$(a)$] For any $i\in{1,2}$, there exists a $j$ such that $\partial_{y_j}f_i(a,b)\neq0$.
\item[$(b)$] If $y_j=p(x,y_{j'})$ is such that $f_i\circ p=0$ in some rectangle of $(a,b_{j'})$, then $\partial_{y_{j'}}(f_{i'}\circ p)(a,b_{j'})\neq0$, where $b=(b_1,b_2)$, $i\neq i'$, and $j\neq j'$.
\end{itemize}
\end{lemma}

According to the proof of Lemma \ref{system lemma}, we can see that by Theorem \ref{prob main1}, $(a)$ yields $y_j=p(x,y_{j'})$, which is the sufficient condition for $(b)$. Additionally, the necessary condition for $(b)$ will follow from the assumption that $f$ is continuously differentiable and $\det J_{f,y}(a,b)\neq0$.

\begin{proof}[Proof of Lemma \ref{system lemma}]
Since the determinant of $J_{f,y}(a,b)$ is not equal to zero, the implicit function theorem implies that there exists a rectangle $R\times V$ containing $(a,b)$, on which the system given by (\ref{system idrv}) has a unique solution for $y$. Furthermore, none of the column vectors of $J_{f,y}(a,b)$ are zero. This means that for any $i$, there exists a $j$ such that $\partial_{y_j}f_{i}(a,b)$ is not equal to zero. Therefore, $(a)$ is proven.

To prove part $(b)$, we simplify the notation by assuming that $i=j=1$. By applying the implicit function theorem to part $(a)$, we obtain the existence of a rectangle $R\times V_2$ containing $(a,b_2)$, an interval $V_1$ of $b_1$, and a function $y_1 = p(x,y_2): R\times V_2 \to V_1$ such that $b_1 = p(a,b_2)$ and $f_1(x,p(x,y_2),y_2) = 0$ for all $(x,y_2) \in R \times V_2$. 
The function $p$ is obtained as the limit of some $\tilde{p}_n$, as proven by Theorem \ref{prob main1}.

First, substitute $y_1=p(x,y_2)$ into $f_2$ and define $h(x,y_2)=f_2(x,p(x,y_2),$ $y_2)$. Then, we have $h(a,b_2)=0$. To prove $(b)$, it suffices to show that $\partial_{y_2}h(a,b_2)\ne0$. Since $\partial_{y_1}f_1(a,b)\ne0$ and $\partial_{y_2}h(a,b_2)=\partial_{y_1}f_2(a,b)\partial_{y_2}p(a,b_2)+\partial_{y_2}f_2(a,b)$ by the chain rule, we obtain 
\begin{equation} \label{sub_jacobian}
	\partial_{y_2}h(a,b_2) 
	= \frac{1}{\partial_{y_1}f_1(a,b)}
	\arraycolsep=4pt\def\arraystretch{1.75}
	\begin{vmatrix}
	\partial_{y_1}f_1(a,b) & 0  \\
	\partial_{y_1} f_2(a,b) & \partial_{y_1}f_2(a,b)\partial_{y_2}p(a,b_2)
			+ \partial_{y_2}f_2(a,b)
	\end{vmatrix}.
\end{equation}

In the determinant of (\ref{sub_jacobian}), after subtracting the first column times $\partial_{y_2}p(a,b_2)$ from the second column, we obtain that the determinant part of (\ref{sub_jacobian}) is equal to
\begin{equation} \label{jacob_column}
	\arraycolsep=4pt\def\arraystretch{1.75}
	\begin{vmatrix}
	\partial_{y_1}f_1(a,b) & -\partial_{y_1}f_1(a,b)\partial_{y_2}p(a,b_2)  \\
	\partial_{y_1} f_2(a,b) & \partial_{y_2}f_2(a,b)
	\end{vmatrix}.
\end{equation}

Since $f_1(x,p(x,y_2),y_2)=0$ holds in $R\times V_2$, we can apply the chain rule to obtain
\begin{equation*}
\begin{split}
	0 &=\partial_{y_2} f_1(a,p(a,b_2),b_2) \\
		&=\partial_{y_1} f_1(a,b)\partial_{y_2}p(a,b_2) + \partial_{y_2} f_1(a,b).
\end{split}	
\end{equation*}
So, (\ref{jacob_column}) is equal to
\begin{equation*}
	\arraycolsep=4pt\def\arraystretch{1.75}
	\begin{vmatrix}
	\partial_{y_1}f_1(a,b) & \partial_{y_2}f_1(a,b)  \\
	\partial_{y_1} f_2(a,b) & \partial_{y_2}f_2(a,b)
	\end{vmatrix}.
\end{equation*}
i.e.,
\begin{equation} \label{last}
	\partial_{y_2}h(a,b_2) 
		=\frac{1}{\partial_{y_1}f_{1}(a,b)} \det J_{f,y}(a,b) \ne 0,
\end{equation}
therefore, the proof is complete.
\end{proof}

In the proof of Lemma \ref{system lemma}, we use the implicit function theorem from the condition of (\ref{last}). This yields a rectangle $R'$ containing $a$, an interval $V_2'$ containing $b_2$, and $y_2=q(x):R'\to V_2'$ such that $h(x,q(x))=f_2(x,p(x,q(x)),q(x))=0$ in $R$, where $q$ is the limit of $\tilde q_n$ by Theorem \ref{prob main1}. We rewrite $R$ and $V_2$ as $R\cap R'$ and $V_2\cap V_2'$, respectively. Then, with $V=V_1\times V_2$, we conclude the main theorem of this section.

\begin{theorem} \label{ift2assumption3}
Using the same notation as in Lemma \ref{system lemma}, if $\det J_{f,y}(a,b) \ne 0$, then there exists a rectangle $R\times V$ containing $(a,b)$,  $y_j=p(x,y_{j'})$ and $y_{j'}=q(x)$ such that $p(x,q(x))$ and $q(x)$ satisfy (\ref{system idrv}) in $R$. Here, $V=V_1\times V_2$, and $p:R\times V_{j'}\to V_j$ and $q:R\to V_{j'}$ are the limits of $\tilde p_n(x,y_{j'})$ and $\tilde q_m(x)$ as $n,m\to\infty$, respectively.
\end{theorem}

The proof of Lemma \ref{system lemma} using Theorem \ref{prob main1} has a natural extension to higher-dimensional versions. Similarly, the proof of Theorems \ref{ift2assumption3} can be extended and applied to systems of equations involving multiple variables for both $x$ and $y$, where the number of equations equals the number of $y$ variables. This approach can be used for implicit functions that are continuously differentiable, without any dimensional restrictions.


\section{Numerical experiment}

We conducted numerical experiments on both an implicit function and a system of two implicit functions. In general, the functions $\tilde g_n$, $\tilde p_n$, and $\tilde q_n$ are step functions that are constant within a partition block. Specifically, if $g$, $p$, and $q$ are polynomials, then they can be identified by $g_n$, $p_n$, and $q_n$ for $n$ greater than or equal to their degree, respectively, as per the implications of Theorem \ref{prob main1}.
Since we are usually interested in analytic implicit functions, we tested these functions using polynomials. We have also presented the pseudo-algorithms \ref{algo1} and \ref{algo2}.
Please note that all entries of a matrix have been rounded off to four decimal places.

First, we provide an example demonstrating the use of Theorem \ref{prob main1}.

%

\def\NoNumber#1{{\def\alglinenumber##1{}\State #1}\addtocounter{ALG@line}{-1}}

\begin{algorithm}[!t]
\caption{for $y=g_n(x)$ of Theorem \ref{prob main1}} \label{algo1}
\begin{algorithmic}[1]
\State\textbf{Preparation} 
\State \quad Let $(a,b)=(a_1,\ldots,a_d,b)$ be such that $f(a,b)=0$.
\State \quad Let $n=$ degree of polynomial.
\State \quad Confirm $\partial_yf(a,b)\ne0$ and fix $\rho=\sign \partial_yf(a,b)$.
\State \quad Take $R\times V$ of $(a,b)$, in which $\rho$ does not change.
\Procedure{implicit}{$f,R,V,(a,b),n$}\Comment{$a\in R$, $b\in V$.}
\State $d\gets\dim R$ \Comment{$d=2$ or $3$.}
\State $R\gets R_1\times\cdots\times R_d$ 
\State $P\gets \frac{1}{2^{nd}}\vert R\vert(\max V)^{\frac{\rho+1}{2}}(\min V)^{\frac{1-\rho}{2}}$
\For{$1\le k\le d$}
\For{$0\le i,j\le2^n-1$}
\State $\delta_{i}(k)\gets \frac{\min R_k + i(\max R_k-\min R_k)}{2^n}$
\State $\Delta_{i}(k)\gets [\delta_{i}(k),\delta_{i+1}(k)]$
\State $v_{i,j}(a_k)\gets \frac{((\delta_{i+1}(k)-a_k)^{j+1}-(\delta_i(k)-a_k)^{j+1})}{j+1}$
\EndFor
\State $A_{n}(a_k)\gets (v_{i,j}(a_k))$
\For{$0\le i_1,\ldots, i_d\le 2^n-1$}
\State $R_{i_1,\ldots,i_d}\gets \Delta_{i_1}(1)\times\cdots\times\Delta_{i_n}(a_d)$
\State $w_{i_1,\ldots,i_d}\gets P-\rho\int_{\min V}^{\max V}\int_{R_{i_1,\ldots,i_d}}\Theta f(x,y)dxdy$
\EndFor
\State $B_n\gets (w_{i_1,\ldots,i_d})$
\EndFor
\For{$1\le k\le d$}
\State $X_k\gets (1,x_k-a_k,(x_k-a_k)^2,\ldots,(x_k-a_k)^{2^n-1})$
\EndFor
\State $(c_{\alpha_1,\ldots,\alpha_d})\gets A_n(a_k)^{-1}\otimes_{k=1}^d B_n$ 
\Comment{$\otimes$ is a tensor contraction.}
   \State \textbf{return} $X_k\otimes_{k=1}^d(c_{\cdots, \alpha_k,\cdots})$ 
   \Comment{$X_1(c_{\alpha_1,\alpha_2}) X_2^T$ if $d=2$.}
\EndProcedure
\end{algorithmic}
\end{algorithm}

\begin{algorithm}[!ht]
\caption{for $y_{j'}=q_m(x)$, $y_j=p_n\circ q_m(x)$ of Theorem \ref{ift2assumption3}}\label{algo2}
\begin{algorithmic}[1]
\State\textbf{Preparation} 
\State \quad Let $(a,b)=(a_1,\ldots,a_d,b_1,b_2)$ be such that $f_1(a,b)=f_2(a,b)=0$.
\State \quad Let $n=$ degree of $p_n$, $m=$ degree of $q_m$.
\State \quad Confirm $\det J_{f,y}(a,b)\ne0$.
\State \quad Take $j$ such that $\partial_{y_j} f_1(a,b)\ne0$ and fix $\rho_j=\sign \partial_{y_j}f_1(a,b)$.
\State \quad Take $R\times V_{j'}$ of $(a,b_j{j'})$, in which $\rho_j$ does not change.
\Procedure{implicit2}{$f_1,f_2,R,V_1,V_2,(a,b),n,m$} \Comment{$a\in R$, $b\in V_1\times V_2$.}
\State $P_n\gets \textsc{implicit}(f_1,R\times V_{j'},V_j,(a,b),n)$ 
\Comment{$(a,b_{j'})\in R\times V_{j'}$, $b_j\in V_j$.}
\State $q_m\gets \textsc{implicit}(f_2\circ P_n,R,V_{j'}(a,b_{j'}),m)$
\Comment{$a\in R$, $b_{j'}\in V_{j'}$.}
\State $p_n\gets P_n\circ q_m$
\Comment{$y_{j'}\gets q_m(x)$.}
   \State \textbf{return} $p_{n}(x)$, $q_m(x)$
   \Comment{$y_j\gets p_{n}(x)$, $y_{j'}\gets q_m(x)$.}
\EndProcedure
\end{algorithmic}
\end{algorithm}

\begin{example} \label{sphere_ex}

We will consider a typical example where the level surfaces of a function $f:\mathbb{R}^2\times\mathbb{R}\to\mathbb{R}$ are spheres. Specifically, we will use the function:
\begin{equation*}
\begin{split}
f(x_1,x_2,y)=x_1^2 + x_2^2 + y^2 - 1,
\end{split}
\end{equation*}
where $f(0,0,1)=0$, as depicted in Figure \ref{fig2}.

\smallskip

\noindent\textbf{Approximation of $y=g(x)$ such that $f(x,y)=0$.}
Let $a=(a_1,a_2)=(0,0)$. Since $\partial_yf(x_1,x_2,y)=2y$, we have $\partial_yf(0,0,1)=2>0$, and therefore $\rho=1$ according to the implicit function theorem. We can select a rectangle $R\times V=[-1/2,1/2)^2\times[0,3/2)$ around the point $(0,0,1)$ (see Figure \ref{fig2}), in which $\rho=1$.

We calculate $y=g_3(x_1,x_2)$ on $[-1/2,1/2)^2$ by applying Theorem \ref{prob main1} with $n=3$ (see Algorithm \ref{algo1}). We obtain the following expressions:
\begin{equation*} 
\begin{split}
	A_3 & \stackrel{\rm{def}}{=} V_3(a_1=0) = V_3(a_2=0) \\
	& = \begin{pmatrix} \left(-1/2+i\Delta\right)^j/j - \left(-1/2+(i-1)\Delta\right)^j/j \end{pmatrix}, \\
	B_3& = 
	\begin{pmatrix} 1/2^{2\cdot 3}\cdot 3/2-\mu(R_{3,(i,j)}) \end{pmatrix}, 
	\end{split}	
\end{equation*}
where $\Delta=1/2^3$ and $1\leq i,j\leq 2^3$.

Using these, we can compute the coefficient matrix of $g_{3}(x_1,x_2)$ at $(0,0)$ as follows:
\begin{equation*} 
\arraycolsep=1.6pt\def\arraystretch{1.2}
A_3^{-1}B_3(A_3^{T})^{-1} =
\begin{matrix}
& \\
\left(
\hspace{0.4em} \vphantom{ \begin{matrix} 12 \\ 12 \\ 12 \\ 12 \\ 12 \\ 12 \\ 12 \\ 12 \end{matrix} } \right .
\end{matrix}
\hspace{-0.8em}
\begin{matrix}
\Scale[1]{1} & \Scale[1]{x_2} & \Scale[1]{x_2^2}  & \Scale[1]{x_2^3} & \Scale[1]{x_2^4}  & \Scale[1]{x_2^5} & \Scale[1]{x_2^6} & \Scale[1]{x_2^7} \\ \arrayrulecolor{gray}\hline \\[-3mm]
\Scale[1]{0.9999}  &  \Scale[1]{0.0}  &  \Scale[1]{-0.5} & \Scale[1]{0.0}  & \Scale[1]{-0.1231}  &  \Scale[1]{0.0} &
\Scale[1]{-0.08}   & \Scale[1]{0.0} \\
\Scale[1]{0.0} &  \Scale[1]{0.0} &  \Scale[1]{0.0} &  \Scale[1]{0.0}  & \Scale[1]{0.0}  &  \Scale[1]{0.0} 
 &  \Scale[1]{0.0} &  \Scale[1]{0.0} \\
\Scale[1]{-0.5}  &  \Scale[1]{0.0} &  \Scale[1]{-0.2504} &  \Scale[1]{0.0} &  \Scale[1]{-0.1808} &  \Scale[1]{0.0} &
\Scale[1]{-0.2237}  & \Scale[1]{0.0} \\
\Scale[1]{0.0} &  \Scale[1]{0.0} &  \Scale[1]{0.0} &  \Scale[1]{0.0}  & \Scale[1]{0.0}  &  \Scale[1]{0.0} 
 &  \Scale[1]{0.0} &  \Scale[1]{0.0} \\
\Scale[1]{-0.1231}  &  \Scale[1]{0.0} &  \Scale[1]{-0.1808} &  \Scale[1]{0.0} &  \Scale[1]{-0.204}  &  \Scale[1]{0.0} &
\Scale[1]{-0.3581}  & \Scale[1]{0.0} \\
\Scale[1]{0.0} &  \Scale[1]{0.0} &  \Scale[1]{0.0} &  \Scale[1]{0.0}  & \Scale[1]{0.0}  &  \Scale[1]{0.0} 
 &  \Scale[1]{0.0} &  \Scale[1]{0.0} \\
\Scale[1]{-0.08}  &  \Scale[1]{0.0} &  \Scale[1]{-0.2237} &  \Scale[1]{0.0} &  \Scale[1]{-0.3581}  &  \Scale[1]{0.0} &
\Scale[1]{-1.3519}  & \Scale[1]{0.0} \\
\Scale[1]{0.0} &  \Scale[1]{0.0} &  \Scale[1]{0.0} &  \Scale[1]{0.0}  & \Scale[1]{0.0}  &  \Scale[1]{0.0} 
 &  \Scale[1]{0.0} & \Scale[1]{0.0}
\end{matrix}
\hspace{-0.2em}
\begin{matrix}
& \\
\left . \vphantom{ \begin{matrix} 12 \\ 12 \\ 12 \\ 12 \\ 12 \\ 12 \\ 12 \\ 12 \end{matrix} } 
\right)
\begin{matrix}
\Scale[1]{1} \\ \Scale[1]{x_1}  \\ \Scale[1]{x_1^2} \\ \Scale[1]{x_1^3} \\ \Scale[1]{x_1^4} \\ \Scale[1]{x_1^5} \\ \Scale[1]{x_1^6} \\ \Scale[1]{x_1^7}
\end{matrix}
\end{matrix},
\end{equation*}
where, for example, $-0.1808$ is the coefficient of $x_1^4x_2^2$ of $g_3$. The surfaces of $g_3$ and  $(1-x_1^2-x_2^2)^{1/2}-g_3(x_1,x_2)$ are plotted in Figures \ref{fig3} and \ref{fig3-1}, respectively, where $(1-x_1^2-x_2^2)^{1/2}$ represents the exact function form of the implicit function.
\begin{figure}[t]
\centerline{\includegraphics[clip, trim=1cm 8cm 1cm 8cm,width=0.8\textwidth]{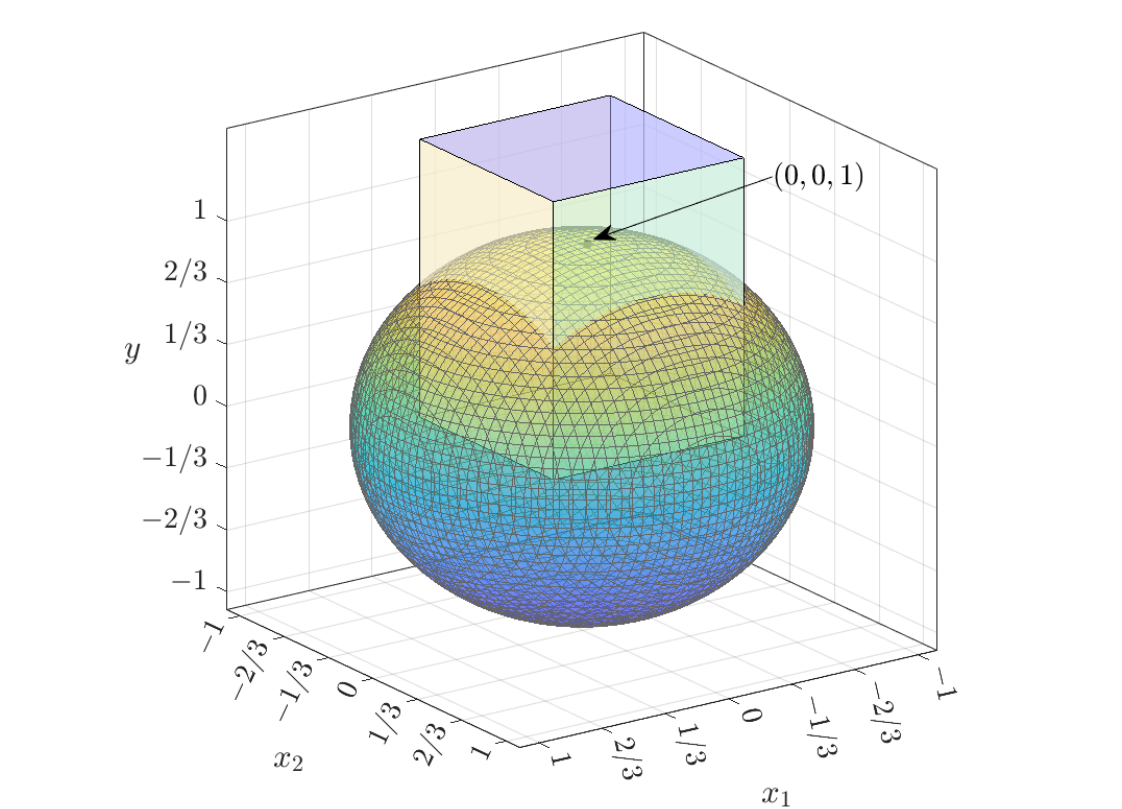}}
\vspace{-0.0cm}
\caption{Surface of $f(x_1,x_2,y)=x_1^2+x_2^2+y^2-1=0$ and $(0,0,1)\in R\times V=[-1/2,1/2)^2\times[0,3/2)$ are illustrated.
}
\label{fig2}
\end{figure}

\begin{figure}[t] 
\centerline{\includegraphics[clip, trim=1cm 8cm 1cm 8cm,width=1\textwidth]{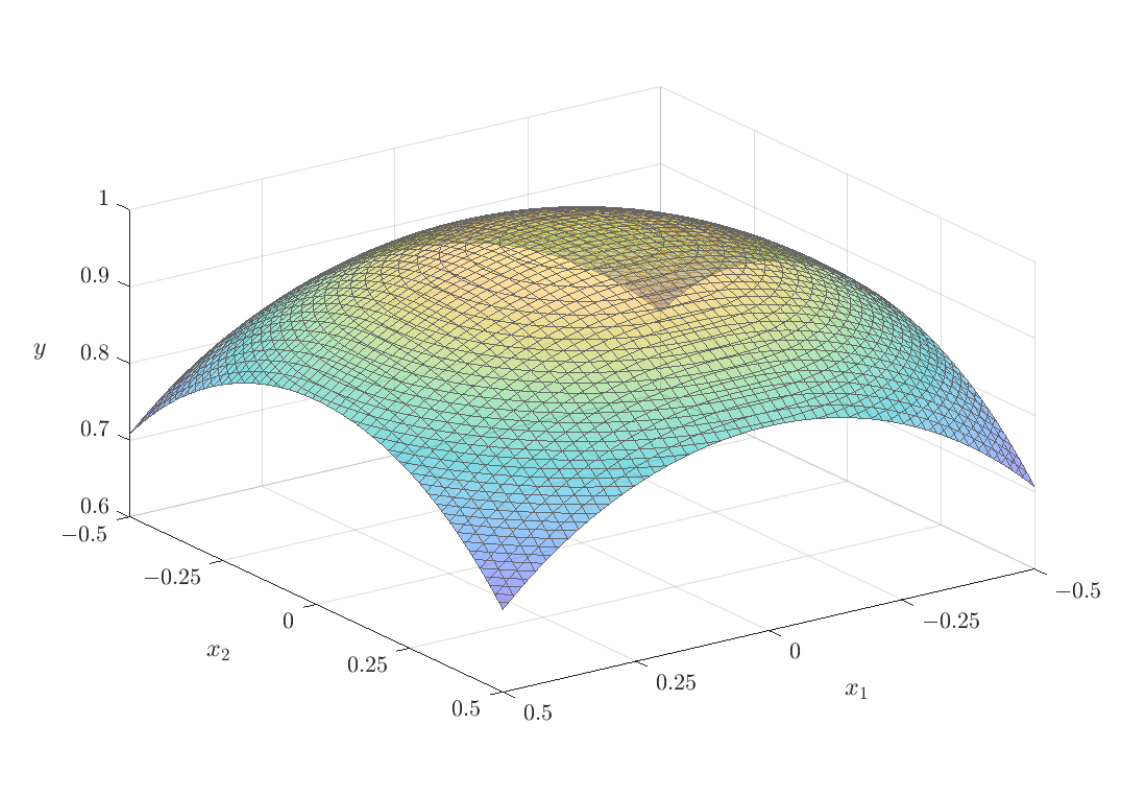} }
\vspace{-0,0cm}    
\caption{Surface $y=g_3(x_1,x_2)$.}
    \label{fig3}%
\end{figure}

\begin{figure}[t] 
\centerline{\includegraphics[clip, trim=1cm 8cm 1cm 8cm,width=1\textwidth]{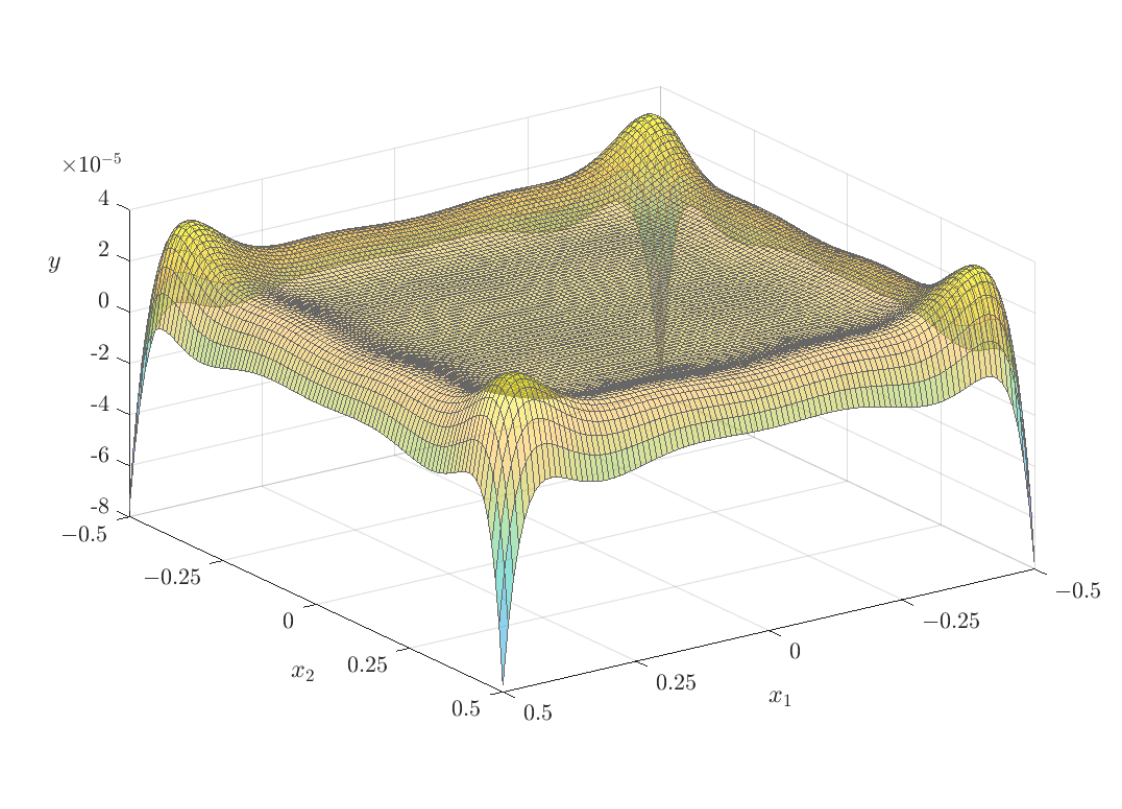} }%
\vspace{-0,5cm}    
\caption{The values of $(1-x_1^2-x_2^2)^{1/2}-g_3(x_1,x_2)$ for accuracy of approximation.}%
    \label{fig3-1}%
\end{figure}
\end{example}

Next, the following example of implicitly defined two equations with three random variables will serve to illustrate Theorem \ref{ift2assumption3}.

\begin{example} \label{nondege}
Suppose that $f_1(x,y_1,y_2)=x + y_1^2 + y_2^3 - 6$ and $f_2(x,y_1,y_2)=x^3y_1 - y_2 - 1$ with $f_1(1,2,1)=f_2(1,2,1)=0$. We aim to approximate the solution functions $y_1=p(x)$ and $y_2=q(x)$ in some interval $R$ containing $x=1$. (The two surfaces of $f_1=0$ and $f_2=0$ are depicted in Figure \ref{fig5-1}.)

\smallskip

\noindent\textbf{Approximation of $y_2=p(x,y_1)$ such that $f_2(x,y_1,y_2)=0$.} 
Let $(a,b_1)=(1,2)$. Since $\det J_{f,y}(1,2,1)=-7\ne0$ ($\det J_{f,y}(x,y_1,y_2)=-2y_1-3x^3y_2^2$), and $\partial_{y_2}f_2(1,2,1)=-1<0$, by the implicit function theorem, there exists a rectangle containing $(1,2,1)$ on which $\rho\stackrel{\rm{def}}{=}\rho_{y_2}=-1$ (To avoid notational ambiguity, we write $\rho$ as $\rho_{y_2}$).
Setting $R\times V_1\times V_2=[1/2,3/2)\times[3/2,5/2)\times[-2,8)$, which is depicted in Figure \ref{fig5-1}, we have $\rho_{y_2}=-1$ in $R\times V_1\times V_2$.

First, for $n=2$ we calculate $y_2=p_{2}(x,y_1)$ such that $f_2=0$ in $R\times V_1$ (refer to Algorithm \ref{algo2}).
Let 
\begin{equation*} 
\begin{aligned}
	A_2(a=1) &= 
	\begin{pmatrix} \left(1/2+i\Delta_x-1\right)^j/j - \left(1/2+(i-1)\Delta_x-1\right)^j/j \end{pmatrix}, \\
	A_2(b_1=2) &= 
	\begin{pmatrix} \left(3/2+i\Delta_{y_1}-2\right)^j/j - \left(3/2+(i-1)\Delta_{y_1}-2\right)^j/j \end{pmatrix}, \\
	B_2 & = 
	\begin{pmatrix} 1/2^{2\cdot2}\cdot (-2)+\mu(R_{2,(i,j)}) \end{pmatrix},
\end{aligned}	
\end{equation*}
where $\Delta_x=\Delta_{y_1}=1/2^2$ and $1\le i,j\le 2^2$.

By applying Theorem \ref{prob main1}, we can calculate the coefficient matrix of $p_2(x,y_1)=\sum_{0\le i,j<4}c_{i,j}(x-1)^{i}(y_1-2)^{j}$ as follows:
\begin{equation} \label{compare10} 
A_2(1)^{-1}B_2(A_2(2)^{T})^{-1} =
\arraycolsep=2pt\def\arraystretch{1.25}
\begin{matrix}
& \\
\left(
\hspace{0.4em} \vphantom{ \begin{matrix} 12 \\ 12 \\ 12 \\ 12  \end{matrix} } \right .
\end{matrix}
\hspace{-0.8em}
\begin{matrix}
\Scale[1]{1} & \Scale[1]{y_1-2} & \Scale[1]{(y_1-2)^2}  & \Scale[1]{(y_1-2)^3}  \\ 
\arrayrulecolor{gray}\hline \\[-5mm]
\Scale[1]{1}   & \Scale[1]{1}  & \Scale[1]{0}  & \Scale[1]{0} \\
\Scale[1]{6}   & \Scale[1]{3}  & \Scale[1]{0}  & \Scale[1]{0} \\
\Scale[1]{6}   & \Scale[1]{3}  & \Scale[1]{0}  & \Scale[1]{0} \\
\Scale[1]{2}   & \Scale[1]{1}  & \Scale[1]{0}  & \Scale[1]{0}
\end{matrix}
\hspace{-0.2em}
\begin{matrix}
& \\
\left . \vphantom{ \begin{matrix} 12 \\ 12 \\ 12 \\ 12 \end{matrix} } 
\right)
\begin{matrix}
\Scale[1]{1} \\ \Scale[1]{x-1}  \\ \Scale[1]{(x-1)^2} \\ \Scale[1]{(x-1)^3} 
\end{matrix}
\end{matrix}
\end{equation}
and $y_2=p_2(x,y_1)$ is drawn in Figure \ref{fig6}. 


\smallskip

\noindent\textbf{Approximation of $y_1=q(x)$ such that $f_1(x_1,y_1,p(x,y_1))=0$.} 
By Theorem \ref{ift2assumption3}, we have $\partial_{y_1}f_1(1,2,p(1,2))\ne0$. Moreover, in $R\times V_1$, the function $f_1(x_1,y_1,p(x,y_1))$ has either $\rho_{y_1}=1$ or $\rho_{y_1}=-1$. Since $y_2=p_2$ is an approximation of $y_2=p$, we can determine the choice of $\rho_{y_1}$ based on the sign of $\partial_{y_1}f_1(1,2,p_2(1,2))$.

With $\rho_{y_1}=1$ determined for $f_1(x,y_1,p_2(x,y_1))$ in $R\times V_1$, we can use a similar method as above to calculate the coefficient array $(c_k)$ of $y_1=q_4(x)$ for $m=4$, which is given by
\begin{equation*}
\begin{split}
(c_k)&=( 2.0021,
-2.5986,
-5.2727,
11.5725,
34.5228,
-61.0522, \\ &\qquad
-144.1976,
164.4454,
362.5225,
-201.1633,
-514.1980, \\ &\qquad
64.6679,
372.6471,
61.3976,
-106.3282,
-38.2053),
\end{split}
\end{equation*}
where
$y_1=q_{4}(x)=\sum_{k=0}^{15}c_k(x-1)^k$ in $[1/2,3/2)$, which satisfies $f_1(x,y_1,$ $p_2(x,y_1))=0$.
(Note that $y_2$ in $f_1(x,y_1,y_2)$ is substituted with $p_2(x,y_1)$ instead of $y_2=p(x,y_1)$.)
The graph of $y_1=q_{4}(x)$ is shown in Figure \ref{fig7}.

%

Finally, we obtain the following approximations for $f_1(x,p(x),q(x))=0$ and $f_2(x,p(x),q(x))=0$ in $[1/2,3/2)$:
\begin{align*}
y_1 &= q_{4}(x) \\
y_2 &= p_{2}(x,q_{4}(x))
\end{align*}
The curve $x\mapsto(q_{4}(x),p_2(x, q_{4}(x)))$ is shown in Figure \ref{fig9}, and the values of $f_k(x,q_{4}(x), p_2(x,q_{4}(x)))$, which approximate 0 for $k=1,2$, are illustrated in Figure \ref{fig8-1}.

\begin{figure}[t]
\centerline{\includegraphics[clip, trim=1cm 7.5cm 1cm 7.5cm, width=0.80\textwidth]{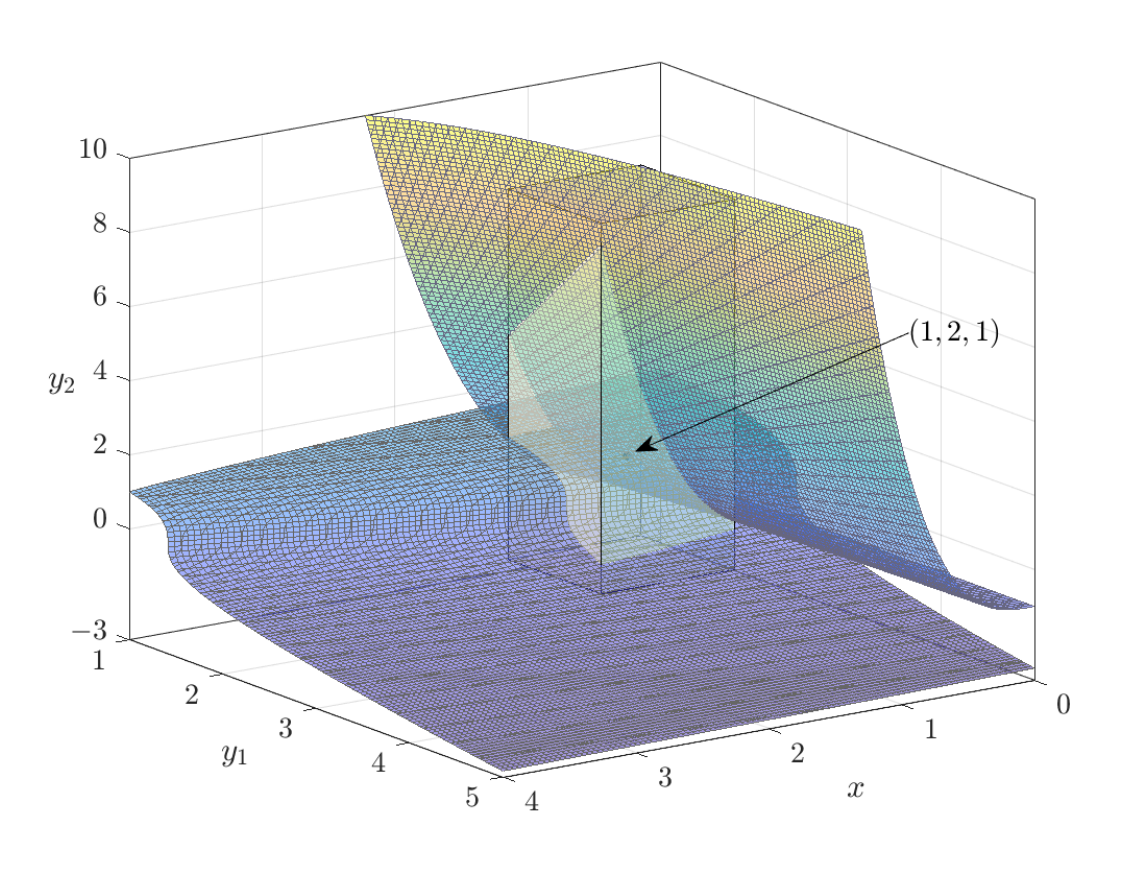}}
\vspace{-0.0cm}
\caption{The surfaces of $f_1(x,y_1,y_2)=x + y_1^2 + y_2^3 - 6=0$ and $f_2(x,y_1,y_2)=x^3y_1 - y_2 - 1=0$ and $R\times V_1\times V_2=[1/2,3/2)\times[3/2,5/2)\times[-2,8)$.}
\label{fig5-1}
\end{figure}

\begin{figure}[H]
\centerline{\includegraphics[clip, trim=1cm 7.5cm 1cm 7.5cm, width=0.80\textwidth]{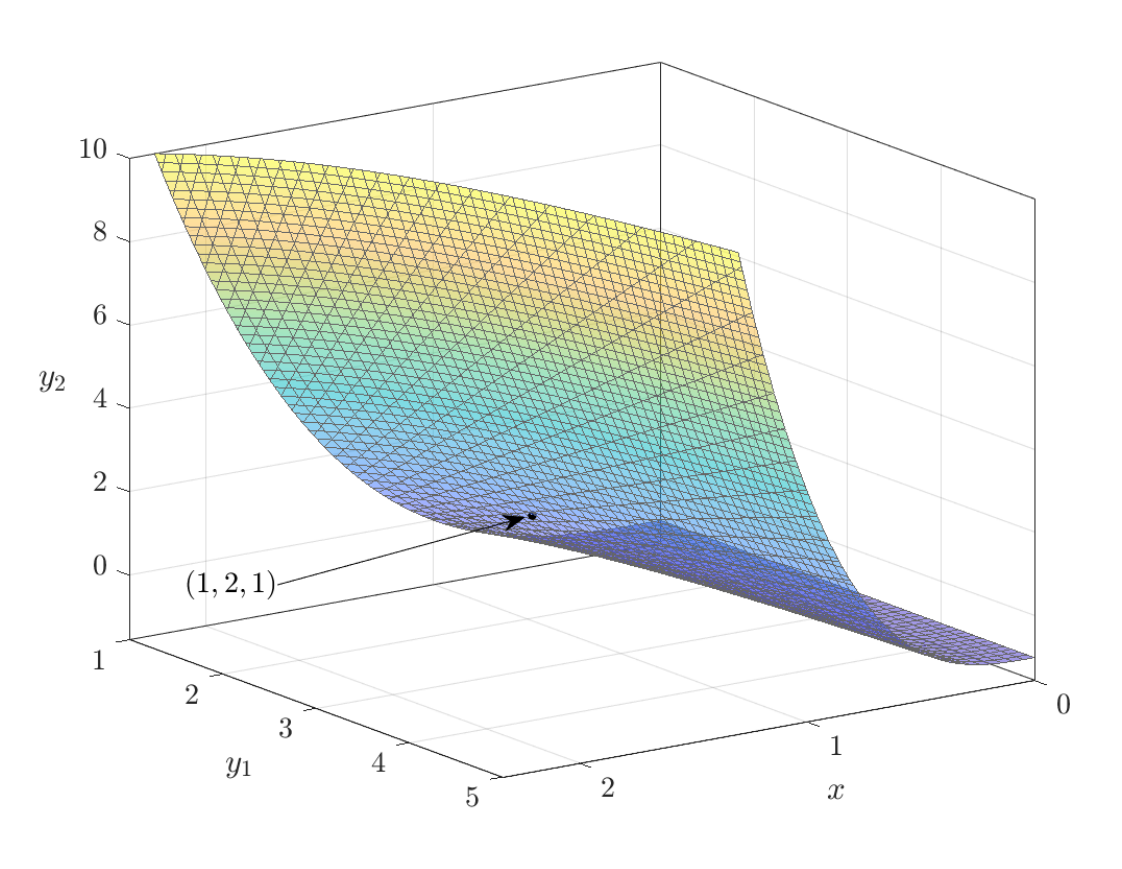}}
\vspace{-0.0cm}
\caption{The surface $y_2=p_2(x,y_1)$.}
\label{fig6}
\end{figure}
\begin{figure}[!ht]
\centerline{\includegraphics[clip, trim=1cm 7.5cm 1cm 7.5cm, width=0.80\textwidth]{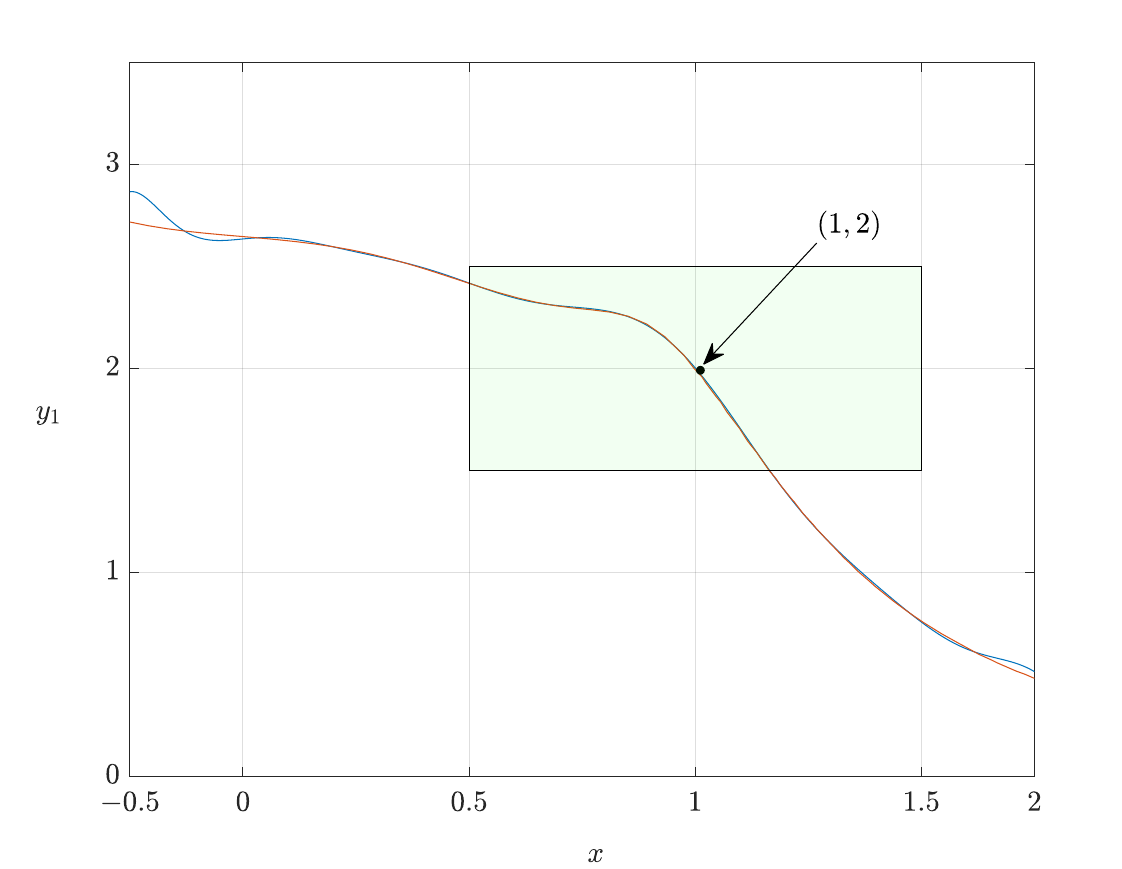}}
\caption{$y_1=q_4(x)$ in blue, which approximates $f_1(x,y_1,p_2(x,y_1))=0$ and a graphical plot of $f_1(x,y_1,p(x,y_1))=0$ in red.}
\label{fig7}
\end{figure}
\begin{figure}[!ht]
\centerline{\includegraphics[clip, trim=1cm 7.5cm 1cm 7.5cm, width=0.85\textwidth]{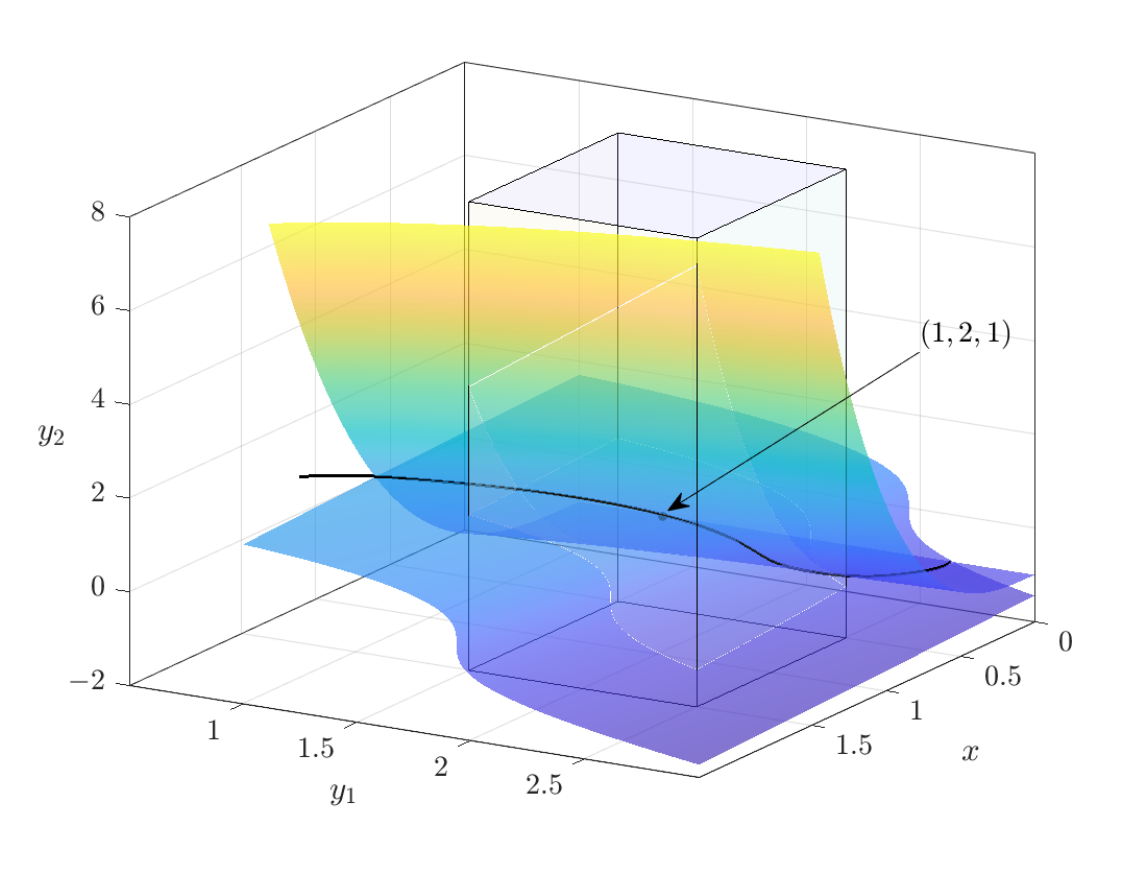}}
\vspace{-0.0cm}
\caption{The curve of $x\mapsto(q_{4}(x),p_2(x,q_{4}(x)))$ is denoted in black.}
\label{fig9}
\end{figure}
\begin{figure}[!ht]
\centerline{\includegraphics[clip, trim=1cm 7.5cm 1cm 7.5cm, width=0.80\textwidth]{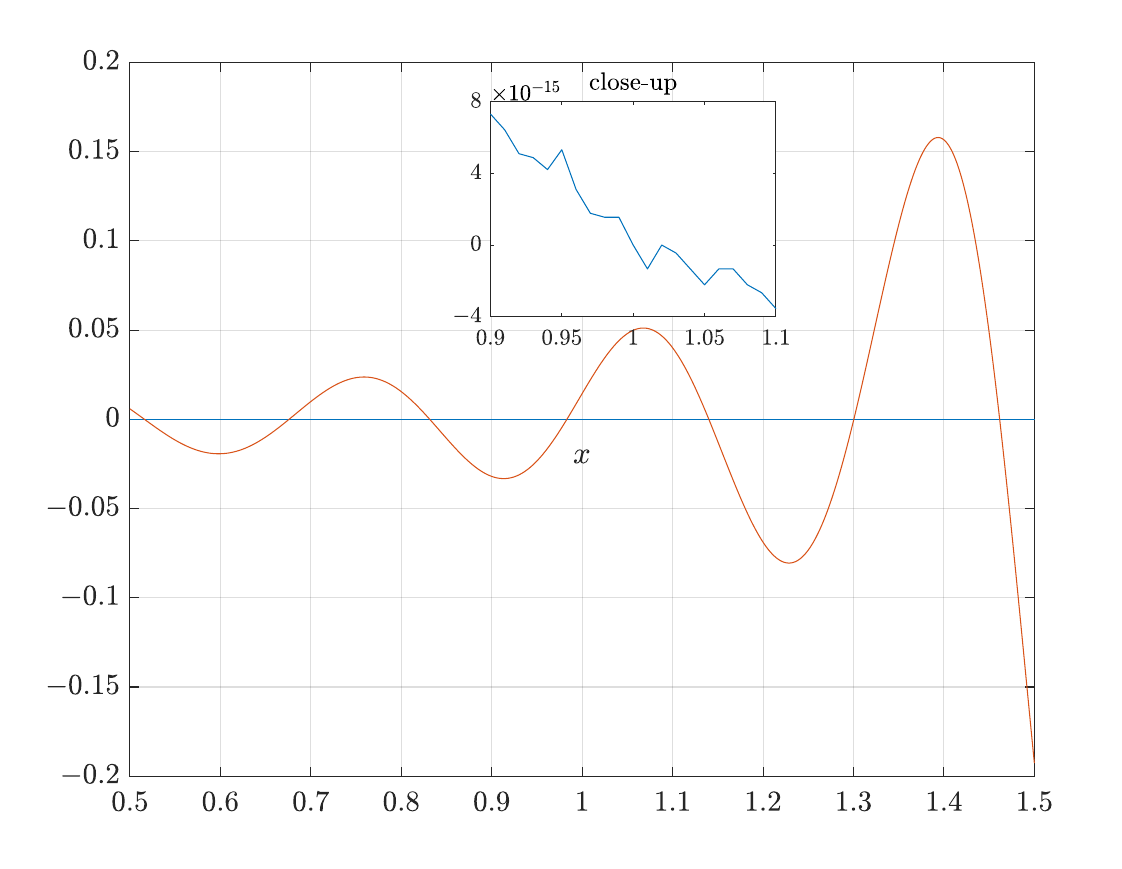}}
\caption{The values of  $f_1(x,q_{4}(x),p_2(x,q_{4}(x)))$ in blue and $f_1(x,q_{4}(x),p_2(x,q_{4}(x)))$ in red on $R=[1/2,3/2)$.}
\label{fig8-1}
\end{figure}
\end{example}

\section*{Conclusion and Further Study}
The implicit function theorem provides local unique existence of a function from an equational expression. To extend the theory of implicit functions, we introduce a polynomial approximation of the implicit function in the $\Gamma$-Ces\`aro sense, based on the integral of Heaviside composition with the condition of nonvanishing Jacobian determinant. The advantage of the proposed method is that it no longer requires higher-order differentiability of an implicit function to obtain higher-order polynomial approximations. Therefore, it is a robust method for implicit functions that are not differentiable to higher-order. We designed two numerical examples to handle an equation with two independent variables and two equations with a single independent variable, but they can be immediately extended to higher dimensions. Additionally, the coefficients of the polynomial are achieved simultaneously over a partition. Finally, the proposed method involves a dyadic decomposition of a domain, and the accuracy of the polynomials relies on domain decomposition. To accelerate the convergence speed, a mathematical strategy for domain decomposition should be studied further.

%

\section*{Declarations}
\begin{itemize}
\item Funding: This work was supported by the National Research Foundation of Korea (No. R1E1A1A03070307).
\item Conflict of interest: The authors declare no competing interests.
\item Ethics approval: No ethical approval was required for this study. 
\item Consent to participate: Agree.
\item Consent for publication: Agree.
\item Availability of data and materials: The datasets are generated and analyzed during the current study.
\item Code availability:  It is available from the corresponding author upon reasonable request. 
\item Authors' contributions: Kyung Soo Rim contributed to mathematical modeling, code development, conceptualization, and literature survey.
\end{itemize}


\end{document}